\documentclass[12pt,reqno]{amsart}
\usepackage{amssymb}
\input amssym.def
\usepackage{amsmath,amsfonts,hyperref,xcolor,mathtools}
\usepackage{amscd}
\usepackage[mathscr]{eucal}
\usepackage{enumitem}
\usepackage{orcidlink}
\allowdisplaybreaks
\numberwithin{equation}{section}

\hypersetup{colorlinks=true,citecolor={purple},linkcolor={teal},urlcolor={violet}}

\theoremstyle{plain}
\newtheorem{theorem}{Theorem}[section]

\newtheorem{lemma}[theorem]{Lemma}

\newtheorem{proposition}[theorem]{Proposition}

\theoremstyle{definition}

\makeatletter
\@namedef{subjclassname@2020}{%
	\textup{2020} Mathematics Subject Classification}
\makeatother

\begin{document}
\title[Linear identities for partition pairs with $5$-cores]{Linear identities for partition pairs with $5$-cores}
\author[Russelle Guadalupe]{Russelle Guadalupe\orcidlink{0009-0001-8974-4502}}
\address{Institute of Mathematics, University of the Philippines Diliman\\
Quezon City 1101, Philippines}
\email{rguadalupe@math.upd.edu.ph}



\begin{abstract}
We prove an infinite family of linear identities for the number $A_5(n)$ of partition pairs of $n$ with $5$-cores by using certain theta function identities involving the Ramanujan's parameter $k(q)$ due to Cooper, and Lee and Park. Consequently, we deduce an infinite family of congruences for $A_5(n)$ using these linear identities.
\end{abstract}

\maketitle

\section{Introduction}\label{sec1}

We denote $f_m := \prod_{n\geq 1}(1-q^{mn})$ for $m\in\mathbb{N}$ and $q\in\mathbb{C}$ with $|q| < 1$ throughout this paper. A partition of $n\in\mathbb{N}$ is a nonincreasing sequence of positive integers $\lambda_1,\ldots,\lambda_m$ such that $\sum_{k=1}^m \lambda_k = n$. The generating function for the number $p(n)$ of partitions of $n$ with $p(0) := 1$ is given by 
\begin{align*}
\sum_{n= 0}^\infty p(n)q^n = \dfrac{1}{f_1}.
\end{align*}
In 1919, Ramanujan \cite{rama} proved the following famous congruences 
\begin{align*}
p(5n+4) &\equiv 0\pmod{5},\\
p(7n+5) &\equiv 0\pmod{7},\\
p(11n+6) &\equiv 0\pmod{11}
\end{align*}
for all $n\geq 0$, which were then extensively refined by Ramanujan \cite{berno}, Watson \cite{wat}, and Atkin \cite{atk} to congruences modulo powers of $5,7$, and $11$ for $p(n)$.  

One may represent a given partition $\lambda = (\lambda_1,\ldots,\lambda_m)$ of $n$ using its Ferris--Young diagram, which is described as follows. We arrange the nodes in $k$ left-aligned rows so that each row has $\lambda_k$ nodes. We then assign to a node at a point $(i,j)$ its hook number, which is the total number of dots directly below and to the right of that node, including the node itself. We call the partition $\lambda$ $t$-core if none of its hook numbers is divisible by $t$ for some $t\in\mathbb{N}$. The generating function for the number $a_t(n)$ of partitions of $n$ that are $t$-cores with $a_t(0) := 1$, due to Garvan, Kim, and Stanton \cite{gks}, is given by 
\begin{align*}
\sum_{n= 0}^\infty a_t(n)q^n = \dfrac{f_t^t}{f_1}.
\end{align*}

Baruah and Berndt \cite{barber} obtained linear identities for $a_3(n)$ and $a_5(n)$ that read
\begin{align}
a_3(4n+1) &= a_3(n),\nonumber\\
a_5(4n+3) &= a_5(2n+1)+2a_5(n)\label{eq11}
\end{align}
for all $n\geq 0$ using Ramanujan's classical modular equations of third \cite[p. 230, Entry 5(i)]{berndt} and fifth degrees \cite[p. 280, Entry 13(iii)]{berndt}, respectively. Kim \cite{kim} found a systematic way via the action of Hecke operators on suitable modular forms of positive weight that yields linear identities for $a_p(n)$ for any prime $p\geq 5$, including (\ref{eq11}). Recently, the author \cite{guad} gave a new elementary proof of (\ref{eq11}) by applying certain identities involving the Ramanujan's parameter
\begin{align*}
k(q) := q\prod_{n=1}^\infty\dfrac{(1-q^{10n-9})(1-q^{10n-8})(1-q^{10n-2})(1-q^{10n-1})}{(1-q^{10n-7})(1-q^{10n-6})(1-q^{10n-4})(1-q^{10n-3})}	
\end{align*}
due to Cooper \cite{coop}, Chern and Tang \cite{chernt}, and Lee and Park \cite{leepark}.

A partition pair of $n$ with $t$-cores is a pair of partitions $(\lambda, \mu)$ such that the sum of all parts of $\lambda$ and $\mu$ is $n$ and both $\lambda$ and $\mu$ are $t$-cores. The generating function for the number $A_t(n)$ of partition pairs of $n$ with $t$-cores is then given by
\begin{align*}
\sum_{n=0}^\infty A_t(n)q^n := \dfrac{f_t^{2t}}{f_1^2}.
\end{align*} 

Baruah and Nath \cite{barnat} used Ramanujan's theta function identities to show that 
\begin{align}\label{eq12}
A_3\left(2^{2k+1}n+\dfrac{5\cdot 2^{2k}-2}{3}\right) = (2^{2k+1}-1)A_3(2n+1),
\end{align}
which subsumes an earlier result of Lin \cite[Theorem 2.6]{lin}. Saikia and Boruah \cite{saibor} found congruences modulo $2$ and $5$ for $A_5(n)$, and Dasappa \cite{dasa} proved an infinite family of congruences modulo powers of $5$ for $A_5(n)$. We refer the interested reader to \cite{barnat,lin,xia,yao,zou} for more congruences for $A_k(n)$ for certain values of $k$.

The goal of this paper is to study arithmetic properties of $A_5(n)$ by relying on elementary $q$-series manipulations. In particular, our main result shows the exact generating function for $A_5(2^kn+2^{k+1}-2)$ for integers $k\geq 1$.

\begin{theorem}\label{thm11}
For integers $k\geq 1$, we have
\begin{align}\label{eq13}
\sum_{n=-1}^\infty A_5(2^kn+2^{k+1}-2)q^n = B_k\dfrac{f_1^4f_5^4}{q}-8B_{k-1}f_2^4f_{10}^4+\dfrac{8^{k+1}-1}{7}\cdot\dfrac{f_5^{10}}{f_1^2}-\dfrac{8^{k+1}-8}{7}\cdot \dfrac{q^2f_{10}^{10}}{f_2^2},
\end{align}
where the sequence $\{B_k\}_{k\geq 0}$ is defined by $B_0=0, B_1=1$, and 
\begin{align*}
B_k = -4B_{k-1}-8B_{k-2}+\dfrac{8^k-1}{7}
\end{align*}
for $k\geq 2$.
\end{theorem}

As a consequence of Theorem \ref{thm11}, we deduce an infinite family of linear identities for $A_5(n)$ analogous to (\ref{eq12}), and an infinite family of congruences for $A_5(n)$.

\begin{theorem}\label{thm12}
For all integers $n\geq 0$ and $k\geq 1$, we have 
\begin{align}\label{eq14}
A_5(2^{k+1}n+3\cdot 2^k-2) =  B_kA_5(4n+4)+\left(\dfrac{8^{k+1}-1}{7}-9B_k\right)A_5(2n+1),
\end{align}
where $\{B_k\}_{k\geq 0}$ is the sequence defined in Theorem \ref{thm11}. Consequently, for all $m\geq 0$ and $n\geq 0$, we have
\begin{align}\label{eq15}
A_5(2^{4m+4}n+3\cdot 2^{4m+3}-2)\equiv 0\pmod{\dfrac{8^{4m+4}-1}{91}}.
\end{align}
\end{theorem}

We organize the remainder of the paper as follows. In Section \ref{sec2}, we present some theta function identities required to establish Theorem \ref{thm11}, which includes the aforementioned identities for $k(q)$ due to Cooper \cite{coop}, and Lee and Park \cite{leepark}, and some $2$-dissection formulas. We then employ these identities to prove Theorem \ref{thm11} in Section \ref{sec3} by finding the generating function for $A_5(2n)$. We finally apply Theorem \ref{thm11} to deduce Theorem \ref{thm12} in Section \ref{sec4}.

\section{Some auxiliary identities}\label{sec2}

We enumerate in this section necessary theta function identities to establish Theorem \ref{thm11}. We begin with the following $2$-dissection formulas and certain identities involving the Ramanujan's parameter $k(q)$.

\begin{lemma}\label{lem21}
We have the identities
\begin{align}
f_1^4 &= \dfrac{f_4^{10}}{f_2^2f_8^4}-4q\dfrac{f_2^2f_8^4}{f_4^2},\label{eq21}\\
\dfrac{f_5}{f_1} &= \dfrac{f_8f_{20}^2}{f_2^2f_{40}}+q\dfrac{f_4^3f_{10}f_{40}}{f_2^3f_8f_{20}}.\label{eq22}
\end{align} 	
\end{lemma}

\begin{proof}
Replacing $q$ with $-q$ in \cite[Lemma 2.3]{hirscr} yields (\ref{eq21}). On the other hand, see \cite[Theorem 2.1]{hirscs} for the proof of (\ref{eq22}).
\end{proof}

\begin{lemma}\label{lem22}
We have the identities
\begin{align}
\dfrac{f_2f_5^5}{qf_1f_{10}^5} &= \dfrac{1}{k(q)}-k(q),\label{eq23}\\
\dfrac{f_2^4f_5^2}{qf_1^2f_{10}^4} &= \dfrac{1}{k(q)}+1-k(q),\label{eq24}\\
\dfrac{f_1^3f_5}{qf_2f_{10}^3} &= \dfrac{1}{k(q)}-4-k(q).\label{eq25}
\end{align}
\end{lemma}

\begin{proof}
See \cite[Theorem 10.4]{coop}.	
\end{proof}

The next identity provides the polynomial relation between $k(q)$ and $k(q^2)$, which was first showed by Lee and Park \cite[Proposition 3.6(1)]{leepark} using modular functions. Recently, the author \cite[Theorem 3.3]{guad} found this relation by applying Lemma \ref{lem22} and an identity involving $k(q)$ and $k(q^2)$ due to Chern and Tang \cite[Theorem 3.3]{chernt}.

\begin{lemma}\label{lem23}
We have the identity
\begin{align*}
X^2-Y+2XY+X^2Y+Y^2=0,	
\end{align*}
where $X:= k(q)$ and $Y := k(q^2)$.
\end{lemma}

The following identity deduces from Lemmas \ref{lem22} and \ref{lem23}, as shown by the author \cite{guad} who used this to illustrate a new proof of (\ref{eq11}). 

\begin{lemma}\label{lem24}
We have the identity
\begin{align*}
\dfrac{f_2^3f_{10}^9}{f_1^3f_4f_5f_{20}^3}-4q^2\dfrac{f_4f_5^2f_{20}^3}{f_1^2} = \dfrac{f_5^5}{f_1}+2q\dfrac{f_{10}^5}{f_2}.
\end{align*}	
\end{lemma}

\begin{proof}
See \cite[Theorem 1]{guad}.	
\end{proof}

We now apply Lemmas \ref{lem22} and \ref{lem23} to derive the following set of theta function identities.

\begin{lemma}\label{lem25}
We have the identity
\begin{align*}
\dfrac{f_4f_{10}^{12}}{f_1^2f_5^2f_{20}^5}+4q^3\dfrac{f_2^3f_5^3f_{20}^5}{f_1^3f_4f_{10}^3} = \dfrac{f_2^3f_5^8}{f_1^4f_{10}^3}-2q\dfrac{f_2^2f_5^3f_{10}^2}{f_1^3}.
\end{align*}
\end{lemma}

\begin{proof}
Dividing both sides of the given identity by $qf_2^2f_5^3f_{10}^2/f_1^3$, it suffices to prove that
\begin{align}\label{eq26}
\dfrac{f_1f_4f_{10}^{10}}{qf_2^2f_5^5f_{20}^5}+4q^2\dfrac{f_2f_{20}^5}{f_4f_{10}^5} = \dfrac{f_2f_5^5}{qf_1f_{10}^5}-2.
\end{align}
Let $X:=k(q)$ and $Y := k(q^2)$.  We replace $q$ with $q^2$ in (\ref{eq23}) so that
\begin{align}\label{eq27}
Z := \dfrac{1}{Y} - Y = \dfrac{f_4f_{10}^5}{q^2f_2f_{20}^5}.
\end{align}
By Lemma \ref{lem23} and (\ref{eq27}), we have
\begin{align*}
1-2X-X^2 = \dfrac{X^2}{Y} + Y = X^2(Z+Y) + Y,
\end{align*}
so that
\begin{align}\label{eq28}
Y = \dfrac{1-2X-X^2-X^2Z}{X^2+1}.
\end{align}	
Combining Lemma \ref{lem23} and (\ref{eq28}) and clearing denominators, we see that
\begin{align}\label{eq29}
X^2Z^2-(1-2X-X^2)(1-X^2)Z+4X(1-X^2)=0.
\end{align}
In view of Lemma \ref{lem22}, (\ref{eq27}), and (\ref{eq29}), we deduce that
\begin{align*}
\dfrac{f_1f_4f_{10}^{10}}{qf_2^2f_5^5f_{20}^5}+4q^2\dfrac{f_2f_{20}^5}{f_4f_{10}^5} &= \dfrac{ZX}{1-X^2}+\dfrac{4}{Z} = \dfrac{XZ^2+4(1-X^2)}{(1-X^2)Z}\\
&=\dfrac{(1-2X-X^2)(1-X^2)Z}{X(1-X^2)Z}=\dfrac{1}{X}-2-X\\
&=\dfrac{f_2f_5^5}{qf_1f_{10}^5}-2,
\end{align*}	
which is exactly (\ref{eq26}). This completes the proof.
\end{proof}

\begin{lemma}\label{lem26}
We have the identity
\begin{align*}
\dfrac{f_1^2f_4^2f_{10}^2}{qf_2^2f_5^2f_{20}^2} - \dfrac{f_2^4f_{20}^2}{f_4^2f_{10}^4} = \dfrac{f_1^3f_5}{qf_2f_{10}^3}. 
\end{align*}	
\end{lemma}

\begin{proof}
Let $X:=k(q)$ and $Y:= k(q^2)$. By Lemma \ref{lem22}, we have 
\begin{align*}
A := \dfrac{1-4X-X^2}{1-X^2} = \dfrac{f_1^4f_{10}^2}{f_2^2f_5^4},\\
B := \dfrac{1-4Y-Y^2}{1-Y^2} = \dfrac{f_2^4f_{20}^2}{f_4^2f_{10}^4}.
\end{align*}
Then from (\ref{eq27}) we have
\begin{align}
\dfrac{4}{1-A} &= \dfrac{1}{X}-X,\label{eq210}\\
Z=\dfrac{4}{1-B} &= \dfrac{1}{Y}-Y.\label{eq211}
\end{align}
Dividing both sides of (\ref{eq29}) by $X^2$ and applying (\ref{eq210}), we get
\begin{align}\label{eq212}
Z^2 - \left(\dfrac{4}{1-A}-2\right)\dfrac{4Z}{1-A}+\dfrac{16}{1-A}=0.
\end{align}
We now substitute (\ref{eq211}) into (\ref{eq212}) and clear denominators, yielding 
\begin{align}\label{eq213}
A^2+4AB+B^2-5A-AB^2=0.
\end{align}
Applying Lemma \ref{lem22}, (\ref{eq210}), and (\ref{eq213}), we arrive at
\begin{align*}
\dfrac{f_1^2f_4^2f_{10}^2}{qf_2^2f_5^2f_{20}^2} - \dfrac{f_2^4f_{20}^2}{f_4^2f_{10}^4} &= \dfrac{A}{B}\left(\dfrac{4}{1-A}+1\right)-B=\dfrac{5A+A^2B-A^2-B^2}{B(1-A)}\\
&=\dfrac{4A}{1-A}=\dfrac{1-X^2}{X}\cdot \dfrac{1-4X-X^2}{1-X^2}\\
&=\dfrac{1}{X}-4-X=\dfrac{f_1^3f_5}{qf_2f_{10}^3}
\end{align*}
as desired.
\end{proof}

\begin{lemma}\label{lem27}
We have the identity	
\begin{align*}
\left(\dfrac{f_4f_{10}^{12}}{f_1^2f_5^2f_{20}^5}-4q^3\dfrac{f_2^3f_5^3f_{20}^5}{f_1^3f_4f_{10}^3}\right)^2 = \dfrac{f_2^4f_5^{12}}{f_1^4f_{10}^4}+4q^2\dfrac{f_2^2f_5^2f_{10}^6}{f_1^2}. 
\end{align*}	
\end{lemma}

\begin{proof}
We know from Lemma \ref{lem25} that
\begin{align}
\left(\dfrac{f_4f_{10}^{12}}{f_1^2f_5^2f_{20}^5}-4q^3\dfrac{f_2^3f_5^3f_{20}^5}{f_1^3f_4f_{10}^3}\right)^2 &= \left(\dfrac{f_4f_{10}^{12}}{f_1^2f_5^2f_{20}^5}+4q^3\dfrac{f_2^3f_5^3f_{20}^5}{f_1^3f_4f_{10}^3}\right)^2 - 16q^3\dfrac{f_2^3f_5f_{10}^9}{f_1^5}\nonumber\\
&=\left(\dfrac{f_2^3f_5^8}{f_1^4f_{10}^3}-2q\dfrac{f_2^2f_5^3f_{10}^2}{f_1^3}\right)^2 - 16q^3\dfrac{f_2^3f_5f_{10}^9}{f_1^5}.\label{eq214}
\end{align}
By expanding and applying (\ref{eq23}) and (\ref{eq25}), we observe that
\begin{align}
\dfrac{f_1^5}{q^3f_2^3f_5f_{10}^9}&\left(\dfrac{f_2^3f_5^8}{f_1^4f_{10}^3}-2q\dfrac{f_2^2f_5^3f_{10}^2}{f_1^3}\right)^2 - 16\nonumber\\
&= \left(\dfrac{f_2f_5^5}{qf_1f_{10}^5}\right)^3-4\left(\dfrac{f_2f_5^5}{qf_1f_{10}^5}\right)^2+4\dfrac{f_2f_5^5}{qf_1f_{10}^5}-16\nonumber\\
&= \left[\left(\dfrac{f_2f_5^5}{qf_1f_{10}^5}\right)^2+4\right]\left(\dfrac{f_2f_5^5}{qf_1f_{10}^5}-4\right)\nonumber\\
&= \left(\dfrac{f_2^4f_5^{10}}{q^2f_1^2f_{10}^{10}}+4\right)\dfrac{f_1^3f_5}{qf_2f_{10}^3}= \dfrac{f_1f_2f_{5}^{11}}{q^3f_{10}^{13}}+4\dfrac{f_1^3f_5}{qf_2f_{10}^3}.\label{eq215}
\end{align}
Multiplying both sides of (\ref{eq215}) by $q^3f_2^3f_5f_{10}^9/f_1^5$ and comparing with (\ref{eq214}), we obtain the desired identity.
\end{proof}

\begin{lemma}\label{lem28}
We have the identity
\begin{align*}
-q\left(\dfrac{f_5^5}{f_1}+2q\dfrac{f_{10}^5}{f_2}\right)^2+f_1^4f_5^4+9q\dfrac{f_5^{10}}{f_1^2}-8q^3\dfrac{f_{10}^{10}}{f_2^2}=\dfrac{f_2^4f_5^{12}}{f_1^4f_{10}^4}+4q^2\dfrac{f_2^2f_5^2f_{10}^6}{f_1^2}. 
\end{align*}
\end{lemma}

\begin{proof}
Expanding the left-hand side of the given identity and then dividing both sides by $q^3f_2^3f_5f_{10}^9/f_1^5$, it remains to prove that
\begin{align}\label{eq216}
\dfrac{f_1^9f_5^3}{q^3f_2^3f_{10}^9}+8\dfrac{f_1^3f_5^9}{q^2f_2^3f_{10}^9}-4\dfrac{f_1^4f_5^4}{qf_2^4f_{10}^4}-12\dfrac{f_1^5f_{10}}{f_2^5f_5}=\dfrac{f_1f_2f_{5}^{11}}{q^3f_{10}^{13}}+4\dfrac{f_1^3f_5}{qf_2f_{10}^3}.
\end{align}
We simplify the left-hand side of (\ref{eq216}) as follows. Letting $X := k(q)$, we infer from Lemma \ref{lem22} that
\begin{align}
\dfrac{f_1^9f_5^3}{q^3f_2^3f_{10}^9} &= \left(\dfrac{f_1^3f_5}{qf_2f_{10}^3}\right)^3 = \left(\dfrac{1}{X}-4-X\right)^3,\label{eq217}\\
\dfrac{f_1^3f_5^9}{q^2f_2^3f_{10}^9} &= \dfrac{f_1^3f_5}{qf_2f_{10}^3}\cdot\dfrac{qf_1^2f_{10}^4}{f_2^4f_5^2}\left(\dfrac{f_2f_5^5}{qf_1f_{10}^5}\right)^2=\left(\dfrac{1}{X}-4-X\right)\dfrac{X}{1+X-X^2}\left(\dfrac{1}{X}-X\right)^2,\label{eq218}\\
\dfrac{f_1^4f_5^4}{qf_2^4f_{10}^4} &= \dfrac{f_1^3f_5}{qf_2f_{10}^3}\cdot\dfrac{qf_1^2f_{10}^4}{f_2^4f_5^2}\cdot \dfrac{f_2f_5^5}{qf_1f_{10}^5}= \left(\dfrac{1}{X}-4-X\right)\dfrac{X}{1+X-X^2}\left(\dfrac{1}{X}-X\right),\label{eq219}\\
\dfrac{f_1^5f_{10}}{f_2^5f_5} &= \dfrac{f_1^3f_5}{qf_2f_{10}^3}\cdot\dfrac{qf_1^2f_{10}^4}{f_2^4f_5^2}=\left(\dfrac{1}{X}-4-X\right)\dfrac{X}{1+X-X^2}.\label{eq220}
\end{align}
Applying (\ref{eq217})--(\ref{eq220}) and Lemma \ref{lem22}, we find that
\begin{align*}
	&\dfrac{f_1^9f_5^3}{q^3f_2^3f_{10}^9}+8\dfrac{f_1^3f_5^9}{q^2f_2^3f_{10}^9}-4\dfrac{f_1^4f_5^4}{qf_2^4f_{10}^4}-12\dfrac{f_1^5f_{10}}{f_2^5f_5}\\
	&=\left(\dfrac{1}{X}-4-X\right)^3+\left(\dfrac{1}{X}-4-X\right)\dfrac{4X}{1+X-X^2}\left[2\left(\dfrac{1}{X}-X\right)^2-\left(\dfrac{1}{X}-X\right)-3\right]\\
	&=\left(\dfrac{1}{X}-4-X\right)^3+\left(\dfrac{1}{X}-4-X\right)\dfrac{4X}{1+X-X^2}\left(\dfrac{2}{X}-3-2X\right)\left(\dfrac{1}{X}+1-X\right)\\
	&=\left(\dfrac{1}{X}-4-X\right)\left[\left(\dfrac{1}{X}-4-X\right)^2+\dfrac{8}{X}-12-8X\right]\\
	&=\left(\dfrac{1}{X}-4-X\right)\left[\left(\dfrac{1}{X}-X\right)^2+4\right]=\dfrac{f_1^3f_5}{qf_2f_{10}^3}\left(\dfrac{f_2^4f_5^{10}}{q^2f_1^2f_{10}^{10}}+4\right)\\
	&= \dfrac{f_1f_2f_{5}^{11}}{q^3f_{10}^{13}}+4\dfrac{f_1^3f_5}{qf_2f_{10}^3},
\end{align*}
which is exactly the right-hand side of (\ref{eq216}). This completes the proof.
\end{proof}

\section{Proof of Theorem \ref{thm11}}\label{sec3}

We establish in this section Theorem \ref{thm11} using the identities derived from Section \ref{sec2}. As an application of these identities, we provide two generating function formulas needed to prove the main result of this paper.

\begin{proposition}\label{prop31}
We have 
\begin{align*}
	\sum_{n=0}^\infty A_5(2n)q^n = f_1^4f_5^4 +9q\dfrac{f_5^{10}}{f_1^2} - 8q^3\dfrac{f_{10}^{10}}{f_2^2}.
\end{align*}
\end{proposition}

\begin{proof}
Using Lemma \ref{lem21}, we expand
\begin{align}\label{eq31}
	\sum_{n=0}^\infty A_5(n)q^n = \dfrac{f_5^{10}}{f_1^2} = \dfrac{f_5^2}{f_1^2}\cdot f_5^8=\left(\dfrac{f_8f_{20}^2}{f_2^2f_{40}}+q\dfrac{f_4^3f_{10}f_{40}}{f_2^3f_8f_{20}}\right)^2\left(\dfrac{f_{20}^{10}}{f_{10}^2f_{40}^4}-4q^5\dfrac{f_{10}^2f_{40}^4}{f_{20}^2}\right)^2.
\end{align}	
We consider the terms in the expansion of (\ref{eq31}) containing $q^{2n}$. In view of Lemmas \ref{lem24}, \ref{lem27}, and \ref{lem28}, we see that
\begin{align*}
	\sum_{n=0}^\infty A_5(2n)q^n &= \dfrac{f_4^2f_{10}^{24}}{f_1^4f_5^4f_{20}^{10}}+q\dfrac{f_2^6f_{10}^{18}}{f_1^6f_4^2f_5^2f_{20}^6}-16q^3\dfrac{f_2^3f_5f_{10}^9}{f_1^5}+16q^5\dfrac{f_4^2f_5^4f_{20}^6}{f_1^4}\\
	&+16q^6\dfrac{f_2^6f_5^6f_{20}^{10}}{f_1^6f_4^2f_{10}^6}\\
	&=\left(\dfrac{f_4f_{10}^{12}}{f_1^2f_5^2f_{20}^5}-4q^3\dfrac{f_2^3f_5^3f_{20}^5}{f_1^3f_4f_{10}^3}\right)^2+q\left(\dfrac{f_2^3f_{10}^9}{f_1^3f_4f_5f_{20}^3}-4q^2\dfrac{f_4f_5^2f_{20}^3}{f_1^2}\right)^2\\
	&=\left(\dfrac{f_4f_{10}^{12}}{f_1^2f_5^2f_{20}^5}-4q^3\dfrac{f_2^3f_5^3f_{20}^5}{f_1^3f_4f_{10}^3}\right)^2+q\left(\dfrac{f_5^5}{f_1}+2q\dfrac{f_{10}^5}{f_2}\right)^2\\
	&=\dfrac{f_2^4f_5^{12}}{f_1^4f_{10}^4}+4q^2\dfrac{f_2^2f_5^2f_{10}^6}{f_1^2}+q\left(\dfrac{f_5^5}{f_1}+2q\dfrac{f_{10}^5}{f_2}\right)^2\\
	&=f_1^4f_5^4 +9q\dfrac{f_5^{10}}{f_1^2} - 8q^3\dfrac{f_{10}^{10}}{f_2^2} 
\end{align*}
as desired.
\end{proof}

\begin{proposition}\label{prop32}
Let $q^{-3}f_1^4f_5^4 := \sum_{n=-3}^\infty c(n+3)q^n$. Then
\begin{align*}
	\sum_{n=-1}^\infty c(2n+3)q^n = -4\dfrac{f_1^4f_5^4}{q}-8f_2^4f_{10}^4.
\end{align*}
\end{proposition}

\begin{proof}
Using (\ref{eq21}), we write
\begin{align}
	\sum_{n=-3}^\infty c(n+3)q^n = \dfrac{f_1^4f_5^4}{q^3}=\dfrac{1}{q^3}\left(\dfrac{f_4^{10}}{f_2^2f_8^4}-4q\dfrac{f_2^2f_8^4}{f_4^2}\right)\left(\dfrac{f_{20}^{10}}{f_{10}^2f_{40}^4}-4q^5\dfrac{f_{10}^2f_{40}^4}{f_{20}^2}\right).\label{eq32}
\end{align}
We extract the terms in the expansion of (\ref{eq32}) involving $q^{2n}$, so that 
\begin{align}\label{eq33}
	\sum_{n=-1}^\infty c(2n+3)q^n = -4\dfrac{f_1^2f_4^4f_{10}^{10}}{qf_2^2f_5^2f_{20}^4}-4q\dfrac{f_2^{10}f_5^2f_{20}^4}{f_1^2f_4^4f_{10}^2}.
\end{align}
We now multiply both sides of Lemma \ref{lem26} by $f_2f_5f_{10}^3/f_1$ and then square both sides of the resulting expression. We obtain
\begin{align}
	\dfrac{f_1^4f_5^4}{q^2} = \left(\dfrac{f_1f_4^2f_{10}^5}{qf_2f_5f_{20}^2}-\dfrac{f_2^5f_5f_{20}^2}{f_1f_4^2f_{10}}\right)^2=\dfrac{f_1^2f_4^4f_{10}^{10}}{q^2f_2^2f_5^2f_{20}^4}-2\dfrac{f_2^4f_{10}^4}{q}+\dfrac{f_2^{10}f_5^2f_{20}^4}{f_1^2f_4^4f_{10}^2}.\label{eq34}
\end{align}
We infer from (\ref{eq33}) and (\ref{eq34}) that
\begin{align*}
	\sum_{n=-1}^\infty c(2n+3)q^n &= -4q\left(\dfrac{f_1^2f_4^4f_{10}^{10}}{q^2f_2^2f_5^2f_{20}^4}+\dfrac{f_2^{10}f_5^2f_{20}^4}{f_1^2f_4^4f_{10}^2}\right)\\
	&= -4q\left(\dfrac{f_1^4f_5^4}{q^2}+2\dfrac{f_2^4f_{10}^4}{q}\right)\\
	&= -4\dfrac{f_1^4f_5^4}{q}-8f_2^4f_{10}^4
\end{align*}
as desired. 
\end{proof}

\begin{proof}[Proof of Theorem \ref{thm11}]
We proceed by induction on $k$. From Proposition \ref{prop31}, we have that
\begin{align}\label{eq35}
	\sum_{n=-1}^\infty A_5(2n+2)q^n = \dfrac{f_1^4f_5^4}{q} +9\dfrac{f_5^{10}}{f_1^2} - 8q^2\dfrac{f_{10}^{10}}{f_2^2},
\end{align}
so (\ref{eq13}) holds for $k=1$. Suppose now that (\ref{eq13}) holds for some $k=m\geq 1$. We divide both sides of (\ref{eq13}) by $q^2$ so that
\begin{align}\label{eq36}
	\sum_{n=-3}^\infty A_5&(2^mn+2^{m+2}-2)q^n\nonumber\\
	&= B_m\dfrac{f_1^4f_5^4}{q^3}-8B_{m-1}\dfrac{f_2^4f_{10}^4}{q^2}+\dfrac{8^{m+1}-1}{7}\cdot\dfrac{f_5^{10}}{q^2f_1^2}-\dfrac{8^{m+1}-8}{7}\cdot \dfrac{f_{10}^{10}}{f_2^2}\nonumber\\
	&= B_m\sum_{n=-3}^\infty c(n+3)q^n-8B_{m-1}\dfrac{f_2^4f_{10}^4}{q^2}+\dfrac{8^{m+1}-1}{7q^2}\sum_{n= 0}^\infty A_5(n)q^n-\dfrac{8^{m+1}-8}{7}\cdot \dfrac{f_{10}^{10}}{f_2^2}.
\end{align}
We extract the terms in the expansion of (\ref{eq36}) involving $q^{2n}$. We deduce from Propositions \ref{prop31} and \ref{prop32} that
\begin{align*}
		\sum_{n=-1}^\infty A_5&(2^{m+1}n+2^{m+2}-2)q^n\\
		&= B_m\sum_{n=-1}^\infty c(2n+3)q^n-8B_{m-1}\dfrac{f_1^4f_5^4}{q}+\dfrac{8^{m+1}-1}{7q}\sum_{n= 0}^\infty A_5(2n)q^n-\dfrac{8^{m+1}-8}{7}\cdot \dfrac{f_5^{10}}{f_1^2}\\
		&= B_m\left(-4\dfrac{f_1^4f_5^4}{q}-8f_2^4f_{10}^4\right)-8B_{m-1}\dfrac{f_1^4f_5^4}{q}+\dfrac{8^{m+1}-1}{7q}\left(f_1^4f_5^4 +9q\dfrac{f_5^{10}}{f_1^2}\right.\\
		&\left.- 8q^3\dfrac{f_{10}^{10}}{f_2^2}\right)-\dfrac{8^{m+1}-8}{7}\cdot \dfrac{f_5^{10}}{f_1^2}\\
		&= \left(-4B_m-8B_{m-1}+\dfrac{8^{m+1}-1}{7}\right)\dfrac{f_1^4f_5^4}{q}-8B_mf_2^4f_{10}^4\\
		&+\left(\dfrac{9(8^{m+1}-1)}{7}-\dfrac{8^{m+1}-8}{7}\right)\dfrac{f_5^{10}}{f_1^2}-\dfrac{8(8^{m+1}-1)}{7}\cdot \dfrac{q^2f_{10}^{10}}{f_2^2}\\
		&= B_{m+1}\dfrac{f_1^4f_5^4}{q}-8B_mf_2^4f_{10}^4+\dfrac{8^{m+2}-1}{7}\cdot\dfrac{f_5^{10}}{f_1^2}-\dfrac{8^{m+2}-8}{7}\cdot \dfrac{q^2f_{10}^{10}}{f_2^2},
\end{align*}
where the last equality follows from the definition of $B_k$. Thus, (\ref{eq13}) holds for $k=m+1$, so it holds for all $k\geq 1$ by induction.
\end{proof}

\section{Proof of Theorem \ref{thm12}}\label{sec4}

As an application of Theorem \ref{thm11}, we prove in this section Theorem \ref{thm12}. We first show the following result, which will be needed to deduce the congruences (\ref{eq15}) stated in the latter theorem. 

\begin{lemma}\label{lem41}
For all $m\geq 0$, we have $B_{4m+3} = \dfrac{8^{4m+4}-1}{91}$.	
\end{lemma}

\begin{proof}
Observe that for $m\geq 0$, 
\begin{align}
	B_{4m+7}+64B_{4m+3}&= B_{4m+7}+4B_{4m+6}+8B_{4m+5}-4(B_{4m+6}+4B_{4m+5}+8B_{4m+4})\nonumber\\
	&+8(B_{4m+5}+4B_{4m+4}+8B_{4m+3})\nonumber\\
	&=\dfrac{8^{4m+7}-1}{7}-4\cdot \dfrac{8^{4m+6}-1}{7}+8\cdot \dfrac{8^{4m+5}-1}{7}=\dfrac{5(8^{4m+6}-1)}{7}.\label{eq41}
\end{align}
By the theory of linear recurrences, we obtain 
\begin{align*}
	B_{4m+3} = A\cdot 8^{4m} + B(-64)^m + C
\end{align*}
for all $m\geq 0$ and some constants $A, B$ and $C$. Using the recurrence relation for $B_k$, we have $B_3=45$, and employing (\ref{eq41}), we compute $B_7=184365$ and $B_{11}=755159085$. Thus, we find $(A,B,C) = (4096/91,0,-1/91)$, leading us to the desired value of $B_{4m+3}$.
\end{proof}

\begin{proof}[Proof of Theorem \ref{thm12}]
We divide both sides of (\ref{eq13}) by $q$ so that
\begin{align}\label{eq42}
	\sum_{n=-2}^\infty A_5(2^kn+3\cdot 2^k-2)q^n = B_k\dfrac{f_1^4f_5^4}{q^2}&-8B_{k-1}\dfrac{f_2^4f_{10}^4}{q}+\dfrac{8^{k+1}-1}{7}\cdot\dfrac{f_5^{10}}{qf_1^2}-\dfrac{8^{k+1}-8}{7}\cdot \dfrac{qf_{10}^{10}}{f_2^2}.
\end{align}
Define $q^{-2}f_1^4f_5^4 :=\sum_{n=-2}^\infty c(n+2)q^n$. We note that the $q$-expansions of $q^{-1}f_2^4f_{10}^4$ and $qf_{10}^{10}/f_2^2$ contain only terms with odd exponents. Thus, by looking at the terms in the expansion of (\ref{eq42}) involving $q^{2n}$, we find that
\begin{align}\label{eq43}
	\sum_{n=-1}^\infty A_5(2^{k+1}n+3\cdot 2^k-2)q^n = B_k\sum_{n=-1}^\infty c(2n+2)q^n+\dfrac{8^{k+1}-1}{7}\sum_{n=0}^\infty A_5(2n+1)q^n.
\end{align}
We next divide both sides of (\ref{eq35}) by $q$ and consider terms in the resulting expansion involving $q^{2n}$. We get
\begin{align}\label{eq44}
	\sum_{n=-1}^\infty A_5(4n+4) = \sum_{n=-1}^\infty c(2n+2)q^n+9\sum_{n=0}^\infty A_5(2n+1)q^n.
\end{align}
Multiplying both sides of (\ref{eq44}) by $B_k$ and subtracting from (\ref{eq43}), we obtain
\begin{align*}
	\sum_{n=-1}^\infty (A_5(2^{k+1}n+3\cdot 2^k-2)-B_kA_5(4n+4))q^n=\left(\dfrac{8^{k+1}-1}{7}-9B_k\right)\sum_{n=0}^\infty A_5(2n+1)q^n.
\end{align*}
Comparing the coefficients of $q^n$ for $n\geq 0$ on both sides of the above expression yields (\ref{eq14}). 

We now set $k=4m+3$ in (\ref{eq14}) and use Lemma \ref{lem41}. Since
\begin{align*}
	\dfrac{8^{4m+4}-1}{7}-9B_{4m+3} = \dfrac{8^{4m+4}-1}{7}-\dfrac{9(8^{4m+4}-1)}{91} = \dfrac{4(8^{4m+4}-1)}{91},
\end{align*}
we finally arrive at
\begin{align*}
	A_5(2^{4m+4}n+3\cdot 2^{4m+3}-2) = \dfrac{8^{4m+4}-1}{91}\left(A_5(4n+4)+4A_5(2n+1)\right)
\end{align*}
for all $m\geq 0$ and $n\geq 0$, which immediately proves (\ref{eq15}).
\end{proof}

\end{document}